%% file: main.tex
\documentclass{article}
\usepackage{amsmath,amsthm,amsfonts,amssymb}
\usepackage[colorlinks]{hyperref}
\usepackage[all]{xy}
\usepackage[pdf]{pstricks}

\title{Spinor sheaves on singular quadrics}
\author{Nicolas Addington}
\date{}

\newcommand\CC{\mathbb C}
\newcommand\EE{\mathcal E}
\newcommand\FF{\mathcal F}
\newcommand\HH{\mathcal H}
\newcommand\OO{\mathcal O}
\newcommand\PP{\mathbb P}
\newcommand\ZZ{\mathbb Z}
\renewcommand\phi\varphi

\newcommand\ev{\text{even}}
\newcommand\odd{\text{odd}}

\newcommand\Cl{{C\ell}}
\DeclareMathOperator\codim{codim}
\DeclareMathOperator\rank{rank}
\DeclareMathOperator\Pin{Pin}
\DeclareMathOperator\Spin{Spin}
\DeclareMathOperator\Hom{Hom}
\DeclareMathOperator\Ext{Ext}
\DeclareMathOperator\tr{tr}

\newtheorem{proposition}{Proposition}[section]
\newtheorem{lemma}[proposition]{Lemma}
\newtheorem*{theorem*}{Theorem}

\begin{document}

\maketitle
\input intro
\input construction
\input dependence
\input dual
\input horrocks
\input stable

\bibliographystyle{plain}
\bibliography{main}

\end{document}

%% file: intro.tex

\begin{abstract}
We define reflexive sheaves on a singular quadric $Q$ that generalize the spinor bundles on smooth quadrics, using matrix factorizations of the equation of $Q$.  We study the first properties of these spinor sheaves, give a Horrocks-type criterion, and show that they are semi-stable, and indeed stable in some cases.
\end{abstract}

\section{Introduction}

On smooth quadric hypersurfaces there exist certain vector bundles called \emph{spinor bundles}, which play a role similar to that of the tautological bundles on Grassmannians.  On a $(2n-1)$-dimensional quadric there is one spinor bundle, of rank $2^{n-1}$.  On a $2n$-dimensional quadric there are two, both of rank $2^{n-1}$.  In low dimensions they coincide with other well-known bundles:  On $Q^1 \cong \PP^1$ it is $\OO(1)$.  On $Q^2 \cong \PP^1 \times \PP^1$ they are $\OO(1,0)$ and $\OO(0,1)$.  On $Q^3$, which is the Lagrangian Grassmannian $LG(2,4)$, it is the tautological quotient bundle.  On $Q^4 \cong G(2,4)$ they are the tautological quotient bundle and the dual of the tautological subbundle.

Ottaviani \cite{ottaviani_spinor} gave geometric and representation-theoretic descriptions of spinor bundles, showed that they are stable, and applied them to moduli spaces of vector bundles on $Q^5$ and $Q^6$.  Later \cite{ottaviani_horrocks} he used them to give a Horrocks-type splitting criterion on smooth quadrics.

Kapranov \cite{kapranov_square} described spinor bundles\footnote{His are the duals of Ottaviani's.  We follow Kapranov's convention, in which the spinor bundles are generated by global sections and the duals have none.} via Clifford algebras and used them to show that the derived category of a smooth quadric is generated by an exceptional collection
\begin{align*}
D^b(Q^{2n-1}) &= \langle \OO(-2n+2),\dotsc,\OO(-1),\OO,S \rangle \\
D^b(Q^{2n}) &= \langle \OO(-2n+1),\dotsc,\OO(-1),\OO,S_+, S_- \rangle
\end{align*}
similar to Be{\u\i}linson's exceptional collection in the derived category of projective space \cite{beilinson}.  Langer \cite{langer} described them via explicit matrix factorizations of $q = x_1 x_2 + \dotsb + x_{2m-1} x_{2m}$ and $q = x_0^2 + x_1 x_2 + \dotsb + x_{2m-1} x_{2m}$ and used them to study the Frobenius morphism on smooth quadrics in characteristic $p>2$.

It is these last two approaches that we generalize to singular quadrics.  Given a linear space $\Lambda$ on a quadric $Q$ defined by a polynomial $q$, we construct a left ideal in the Clifford algebra of $q$, and from this a matrix factorization of $q$, and from this two sheaves $S$ and $T$ on $Q$, which we call spinor sheaves.  They fail to be vector bundles where $\Lambda$ meets the singular locus $Q_\text{sing}$.  In contrast to the one or two rigid bundles we had when $Q$ was smooth, now we have one or two families of reflexive sheaves for each Grassmannian of linear spaces on $Q_\text{sing}$.

We treat smooth and singular quadrics uniformly, but even on smooth quadrics our description of spinor bundles has advantages over Ottaviani's, with which it is difficult to do homological algebra, Kapranov's, with which it is difficult to do geometry, and Langer's, with which it is difficult to vary the quadric in a family.

In \S\ref{construction} we give the details of the construction.  In \S\ref{dependence} we describe how $S$ and $T$ vary with $\Lambda$.  In \S\ref{dual} we study their dual sheaves. In \S\ref{horrocks} we describe how they restrict to a hyperplane section of $Q$ and pull back to a cone on $Q$, and we prove a Horrocks-type criterion.  In \S\ref{stable} we show that they are stable when $\Lambda$ is maximal and properly semi-stable otherwise. 

I was motivated to define spinor sheaves by studying a moduli problem on the complete intersection of four quadrics \cite{addington_quadrics}.  This work was supported in part by the National Science Foundation under grants nos.\ 
DMS-0354112, 
DMS-0556042, 
and DMS-0838210. 

%% file: construction.tex

\section{The Construction}\label{construction}

Let $V$ be a complex vector space equipped with a quadratic form $q$ of rank at least 2, so the corresponding quadric hypersuface $Q \subset \PP V$ is reduced.  Let $b$ be the symmetric bilinear form associated to $q$.  Let
\begin{align*}
\Cl &= \mathcal T(V)/\langle v^2 = q(v) \rangle \\
&= \mathcal T(V)/\langle vv' + v'v = 2b(v,v') \rangle
\end{align*}
be the Clifford algebra of $q$.  If $\{ v_1, \dotsc, v_n \}$ is a basis of $V$ then
\[ \{ v_{i_1} \dotsm v_{i_k} : 1 \le i_1 < \dotsb < i_k \le n \} \]
is a basis of $\Cl$.

Let $W \subset V$ be an isotropic subspace, that is, one with $q|_W = 0$, or equivalently $\PP W \subset Q$.  Choose a basis $w_1, \dotsc, w_m$, and let $I$ be the left ideal $I = \Cl \cdot w_1 \dotsm w_m$.  Since $W$ is isotropic, choosing a different basis just rescales the generator $w_1 \dotsm w_m$ by the determinant of the change-of-basis matrix, so $I$ is independent of this choice.  Since $\Cl$ is $\ZZ/2$-graded, we can write $I = I_\ev \oplus I_\odd$.  We will always consider $I$ in the category of graded left $\Cl$-modules and maps that respect the grading.
\vspace\baselineskip

A matrix factorization of $q$ is a pair of $N \times N$ matrices $\phi$ and $\psi$, of polynomials in general but of linear forms in our case, such that
\[ \phi \psi = \psi \phi = q \cdot 1_{N \times N}, \]
where $1_{N \times N}$ is the identity matrix.  For example,
\[ \left[\begin{smallmatrix}
x_3 & x_4 & x_5 & 0 \\
-x_1 & x_0 & 0 & x_5 \\
-x_2 & 0 & x_0 & -x_4 \\
0 & -x_2 & x_1 & x_3
\end{smallmatrix}\right]
\left[\begin{smallmatrix}
x_0 & -x_4 & -x_5 & 0 \\
x_1 & x_3 & 0 & -x_5 \\
x_2 & 0 & x_3 & x_4 \\
0 & x_2 & -x_1 & x_0
\end{smallmatrix}\right]
 = (x_0 x_3 + x_1 x_4 + x_2 x_5) \cdot 1_{4 \times 4}. \]
An $N \times N$ matrix of linear forms is the same as a map $\OO_{\PP V}^N(-1) \to \OO_{\PP V}^N$, and it will be more natural to work with the latter.

From the module $I$, define a map of vector bundles
\[ \begin{array}{rclcrcl}
\OO_{\PP V}(-1) & \otimes & I_\ev & \xrightarrow\phi & \OO_{\PP V} & \otimes & I_\odd \\[2pt]
v & \otimes & \xi & \mapsto & 1 & \otimes & v\xi.
\end{array} \]
Here we are regarding $\OO_{\PP V}(-1)$ as the tautological line bundle, that is, as a subbundle of $\OO_{\PP V} \otimes V$.  Define $\psi: \OO_{\PP V}(-1) \otimes I_\odd \to \OO_{\PP V} \otimes I_\ev$ similarly.  Then the compositions
\[ \begin{array}{rclcrclcrclc}
\OO_{\PP V} & \otimes & I_\ev & \xrightarrow\phi &
\OO_{\PP V}(1) & \otimes & I_\odd & \xrightarrow\psi &
\OO_{\PP V}(2) & \otimes & I_\ev \\[1ex]
\OO_{\PP V} & \otimes & I_\odd & \xrightarrow\psi &
\OO_{\PP V}(1) & \otimes & I_\ev & \xrightarrow\phi &
\OO_{\PP V}(2) & \otimes & I_\odd
\end{array} \]
are just multiplication by $q$, so we have a matrix factorization of $q$.

This link between $\Cl$-modules and matrix factorizations was first studied by Buchweitz, Eisenbud, and Herzog \cite{beh} and has been rediscovered more than once \cite[\S7.4]{kapustin_li} \cite{bertin}.  It is interesting to note the resemblance to the Thom class in K-theory \cite[App.~C]{lawson_michelsohn}.  \vspace\baselineskip

Finally, let $S = \operatorname{coker} \phi$ and $T = \operatorname{coker} \psi$.  These are supported on $Q$, for $\phi$ and $\psi$ are isomorphisms where $q \ne 0$.  Since $\phi \circ \psi$ and $\psi \circ \phi$ are injective, $\phi$ and $\psi$ are injective, so $S$ and $T$ have resolutions on $\PP V$
\begin{equation}\label{main_res}
\begin{split}
&0 \to \OO_{\PP V}^N(-1) \xrightarrow\phi \OO_{\PP V}^N \to S \to 0 \\
&0 \to \OO_{\PP V}^N(-1) \xrightarrow\psi \OO_{\PP V}^N \to T \to 0,
\end{split}
\end{equation}
where $N = \dim I_\ev = \dim I_\odd = 2^{\codim W - 1}$.  From these it is easy to compute their cohomology.

We ask how far $S$ and $T$ are from being vector bundles.  Let $K \subset V$ be the kernel of $q$ and recall that the singular locus of $Q$ is $\PP K$.
\begin{proposition}\label{not_a_vb}
The restriction of $S$ to $\PP K \cap \PP W$ is trivial of rank $2^{\codim W - 1}$.  If $\codim W > 1$ then elsewhere on $Q$, $S$ is locally free of rank $2^{\codim W - 2}$.  The same is true of $T$.
\end{proposition}
\begin{proof}
For each $v \in V$, we want to know the rank of the linear map
\[ I_\ev \xrightarrow{v\cdot} I_\odd \]
given by left multiplication by $v$.

If $v \in W \cap K$ then the map is zero: any $v \in K$ commutes with elements of $\Cl_\ev$ and anti-commutes with elements of $\Cl_\odd$, and any $v \in W$ annihilates the generator $w_1 \dotsm w_m$ of $I$, so any $v \in W \cap K$ annihiliates $I$.

If $v \notin W$, choose a basis of $V$ starting with $v$ and ending with $w_1, \dotsc, w_m$.  Then any element of $I$ can be written uniquely as $(v\xi + \eta) w_1 \dotsm w_m$, where $v$ and $w_1, \dotsc, w_m$ do not appear in $\xi$ and $\eta$.  Now
\[ v \cdot (v\xi + \eta) w_1 \dotsm w_m = v \eta w_1 \dotsm w_m, \]
which is zero if and only if $\eta = 0$.  Thus the rank of the map is $2^{\codim W - 2}$, so $S|_{Q \setminus \PP W}$ is a vector bundle of rank $2^{\codim W - 2}$.  But in \S\ref{dual} we will see that $\EE xt^{>0}_Q(S,\OO_Q) = 0$, which implies that $S$ is a vector bundle on the whole smooth locus of $Q$.
\end{proof}

\noindent If $\codim W = 1$ then $Q = \PP W \cup \PP W'$, where $W'$ is the other maximal isotropic subspace, and it is easy to check that $S = \OO_{\PP W}$ and $T = \OO_{\PP W'}$ when $\dim V$ is odd and vice versa when $\dim V$ is even.

%% file: dependence.tex

\section{Dependence on the Isotropic Subspace}\label{dependence}

In this section we first summarize how $S$ and $T$ vary with $W$ and how they are related for $W$s of various dimensions, then prove the corresponding statements about graded $\Cl$-modules, and finally show that the functor from modules to sheaves is fully faithful.  \vspace\baselineskip

If $q$ is non-degenerate then $S$ and $T$ are rigid, so varying $W$ continuously leaves them unchanged.  If $\dim W < \frac12 \dim V$, so $W$ belongs to a connected family, then $S \cong T$.  If $\dim W = \frac12 \dim V$, so $W$ belongs to one of two families, then $S \not\cong T$, and switching $W$ to the other family interchanges $S$ and $T$.  If $W$ is maximal then when $\dim V$ is odd, $S \cong T$ is the classical spinor bundle, and when $\dim V$ is even, $S$ and $T$ are the two spinor bundles.  If $W'$ is codimension 1 in a maximal $W$ then $S' \cong T' \cong S \oplus T$, if codimension 2 then $S' \cong T' \cong (S \oplus T)^{\oplus 2}$, and in general $S' \cong T' \cong (S \oplus T)^{\oplus 2^{\dim W/W' - 1}}$.  In short then, on smooth quadrics, maximal isotropic subspaces give the classical spinor bundles, and non-maximal ones give direct sums of them.

If $q$ is degenerate then $S$ and $T$ are not rigid in general, since by Proposition \ref{not_a_vb} we can recover $W \cap K$ from them; but varying $W$ continuously while keeping $W \cap K$ fixed leaves them unchanged.  Let $\pi: V \to V/K$ be the projection, and recall that $q$ descends to a non-degenerate form on $V/K$.  If $\dim \pi(W) < \frac12 \dim V/K$ then $S \cong T$.  If $\dim \pi(W) = \frac12 \dim V/K$ then $S \not\cong T$, and switching $\pi(W)$ to the other family (still keeping $W \cap K$ fixed) interchanges them.  For example, consider a line $\PP W$ on a rank 2 quadric surface:
\begin{center}
\psset{unit=5pt}
\begin{pspicture}(14,10)
\pscustom{
  \newpath
  \moveto(10,2)
  \rlineto(4,-1)
  \rlineto(-6,6)
  \rlineto(-8,2)
  \rlineto(3,-3)
  \rlineto(8,-2)
  \rlineto(2,4)
  \rlineto(-8,2)
  \rlineto(-1,-2)
  \moveto(3,6)
  \rlineto(-2,-4)
  \rlineto(8,-2)
  \rlineto(2,4)
}
\pscustom[linestyle=dashed,dash=2pt 2pt]{
  \newpath
  \moveto(4,8)
  \rlineto(-1,-2)
  \rlineto(3,-3)
  \rlineto(8,-2)
}
\pscustom[linewidth=.3pt]{
  \newpath
  \moveto(5,1)
  \rlineto(6,7.5)
}
\end{pspicture}
\end{center}
If $\PP W$ is not the cone line, it lies on one plane or the other and meets the cone line in a point, and these data determine the isomorphism classes of $S$ and $T$: varying the line while keeping the point fixed leaves $S$ and $T$ unchanged, switching to the other plane interchanges them, and varying the point deforms them.

Our earlier comments on direct sums are generalized as follows: if $W'$ is codimension 1 in $W$ then there are exact sequences
\begin{align}\label{jordan_hoelder}
0 \to S \to S' \to T \to 0 & & 0 \to T \to T' \to S \to 0
\end{align}
which split if and only if $W' \cap K = W \cap K$.  That is, if $W \cap K$ shrinks we get interesting extensions, but if $\pi(W)$ shrinks we just get direct sums.  So while on smooth quadrics only maximal $W$s were interesting, on singular quadrics non-maximal $W$s may be interesting, but only if $\pi(W)$ is maximal.  \vspace\baselineskip

To prove all this, we introduce the following group action.  Let $G$ be the subgroup of the group of units $\Cl^\times$ generated by the unit vectors, that is, by those $u \in V$ with $q(u) = 1$, and let $G$ act on $V$ by reflections:
\[ u \cdot v := v - 2b(v,u) u = -uvu^{-1}. \]
This preserves $q$, for
\[ q(-uvu^{-1}) = (-uvu^{-1})^2 = uv^2u^{-1} = u q(v) u^{-1} = q(v). \]
The spinor sheaves $S$ and $T$ are equivariant for the action of $G_\ev := G \cap \Cl_\ev$ on $Q$, as we see from the commutative diagram
\begin{equation}\label{equivariant}
\xymatrix@R=.75\baselineskip@C=\baselineskip{
& \OO_{\PP V}(-1) \otimes I_\ev \ar[rr]^-\phi \ar[dddd] & & \OO_{\PP V} \otimes I_\odd \ar[dddd] \\
v \otimes \xi \ar@{|->}[dd] & & & & 1 \otimes \xi \ar@{|->}[dd] \\ \\
gvg^{-1} \otimes g\xi & & & & 1 \otimes g\xi \\
& \OO_{\PP V}(-1) \otimes I_\ev \ar[rr]^-\phi & & \OO_{\PP V} \otimes I_\odd.
}
\end{equation}

If $q$ is non-degenerate then $G$ is $\Pin(V,q)$, the central extension of the orthogonal group $O(V,q)$ by $\ZZ/2$, and $G_\ev$ is $\Spin(V,q)$ \cite[\S I.2]{lawson_michelsohn}.  When $m < \frac12 \dim V$, $G_\ev$ acts transitively on the variety of $m$-dimensional isotropic subspaces, and when $m = \frac12 \dim V$, $G_\ev$ acts transitively on each of its connected components and $G_\odd$ interchanges them.

If $q$ is degenerate, the natural map $G \to O(V,q)$ is not surjective, since $O(V,q)$ acts transitively on $K$ while $G_\ev$ acts as the identity on $K$ and $G_\odd$ acts as $-1$.  If $U \subset V$ is a subspace complementary to $K$ then $q|_U$ is non-degenerate and $G$ contains $\Pin(U,q|_U)$.  From this it is not hard to see that $G$ can take $W$ to $W'$ if $W \cap K = W' \cap K$ and that $G_\ev$ can if in addition $\pi(W)$ and $\pi(W')$ lie in the same family.

Now if $g \in G$ takes an isotropic subspace $W \subset V$ to another one $W' = gWg^{-1}$, then right multiplication by $g^{-1}$ takes $I$ to $I'$:
\[ \Cl \cdot w_1 \dotsm w_m \cdot g^{-1} = \Cl \cdot (\pm g w_1 g^{-1}) \dotsm (\pm g w_m g^{-1}). \]
So if $g \in G_\ev$ then $I \cong I'$, and if $g \in G_\odd$ then $I \cong I'[1]$, where $I'[1]$ means $I'$ with the odd and even pieces interchanged.  Thus we have proved:
\begin{proposition} \ 
\begin{itemize}
\item If $\dim \pi(W) < \frac12 \dim V/K$ then $I \cong I[1]$.
\item Suppose that $W \cap K = W' \cap K$.  If $\pi(W)$ and $\pi(W')$ lie in the same family then $I \cong I'$.  If they lie in opposite families then $I \cong I'[1]$.
\end{itemize}
\end{proposition}

\noindent Inversely,
\begin{proposition} \ 
\begin{itemize}
\item If $\dim \pi(W) = \frac12 \dim V/K$ then $I \not\cong I[1]$.
\item If $W \cap K \ne W' \cap K$ then $I$ is not isomorphic to $I'$ or $I'[1]$.
\end{itemize}
\end{proposition}
\begin{proof}
For the first statement, let $\dim V/K = 2k$.  Then there is a basis $v_1, \dotsc, v_n$ of $V$ in which
\[ q = x_1 x_{k+1} + \dotsb + x_k x_{2k} \]
and $W = \operatorname{span}(v_{k+1}, \dotsc, v_{2k}, v_{2k+1}, \dotsc, v_{2k+l})$, where $l = \dim(W \cap K)$.  Observe that $\xi := v_1 \dotsm v_k$ annihilates every element of the associated basis of $I$ except $v_{k+1} \dotsm v_{2k} v_{2k+1} \dotsm v_{2k+l}$.  Thus if $\dim W$ is even then $\xi$ annihilates $I_\odd$ but not $I_\ev$, and vice versa if $\dim W$ is odd.  For the second statement, we saw in the proof of Proposition \ref{not_a_vb} that $W \cap K = V \cap \operatorname{Ann}I$, where the latter intersection takes place in $\Cl$.
\end{proof}

\begin{proposition}
Suppose that $W' \subset W$ is codimension 1.  Then for any $w \in W \setminus W'$ the sequence
\[ 0 \to I \to I' \xrightarrow{\cdot w} I[1] \to 0 \]
is exact.  It is split if and only if $W \cap K = W' \cap K$.
\end{proposition}
\begin{proof}
To see that the sequence is exact, choose a basis $w_1, \dotsc, w_m$ of $W'$, and extend this to a basis of $V$ ending with $w, w_1, \dotsc, w_m$.  Then just as we argued in the proof of Proposition \ref{not_a_vb}, an element of $I'$ can be written as $(\xi + \eta w)w_1 \dotsb w_m$, where $w$ and $w_1, \dotsc, w_m$ do not appear in $\xi$ and $\eta$, and if this times $w$ equals zero then $\xi = 0$.

If $W \cap K = W' \cap K$ then $\pi(w) \notin \pi(W')$, so there is a $v \perp W'$ with $b(v,w) = \frac12$.  Since $v \perp W'$, we have $I \cdot v \subset I'$.  Since $b(v,w) = \frac12$, the map $I[1] \xrightarrow{\cdot v} I'$ splits $I' \xrightarrow{\cdot w} I[1]$:
\[ \xi w w_1 \dotsm w_m \cdot v w = \xi w w_1 \dotsm w_m (1 - wv) = \xi w w_1 \dotsm w_m. \]
Inversely, if $W \cap K \ne W' \cap K$ then $I'$ and $I \oplus I[1]$ have different annihilators, hence are not isomorphic.
\end{proof}
\ 

To see that what we have proved about modules implies what we have claimed about spinor sheaves, we study the functor that sends a graded $\Cl$-module $I$ to a sheaf $S$ on $Q$.  It is indeed a functor, for a homogeneous map $f: I \to I'$ induces a commutative diagram
\[ \xymatrix{
\OO_{\PP V}(-1) \otimes I_\ev \ar[r]^-\phi \ar[d]^{1 \otimes f_\ev}
& \OO_{\PP V} \otimes I_\odd \ar[d]^{1 \otimes f_\odd} \\
\OO_{\PP V}(-1) \otimes I'_\ev \ar[r]^-{\phi'} & \OO_{\PP V} \otimes I'_\odd
} \]
and hence a map on cokernels.  The functor is exact.  Kuznetsov \cite[Prop.~4.9]{kuznetsov_quadrics}, working more generally with quadric fibrations, showed that it is fully faithful.  We give a different proof:
\begin{proposition}\label{fully_faithful}
The natural map $\Hom_\Cl(I,I') \to \Hom_Q(S,S')$ is an isomorphism.
\end{proposition}
\begin{proof}
The inverse is essentially the map $\Hom_Q(S,S') \to \Hom_\CC(I_\odd, I'_\odd)$, where the second object is vector space homomorphisms, given by taking global sections.  This is injective because $S$ is generated by global sections.  The composition $\Hom_\Cl(I,I') \to \Hom_\CC(I_\odd, I'_\odd)$ sends $f$ to $f_\odd$.  This too is injective: a map $f$ of graded $\Cl$-modules is determined by $f_\odd$, for if $m \in I_\ev$ and $v \in V$ has $q(v) = 1$ then $f(m) = v^2 f(m) = v f(vm)$.

It remains to check that a linear map $I_\odd \to I'_\odd$ induced by a sheaf map $S \to S'$ is induced by a module map $I \to I'$.  Applying $\Hom_Q(-,S')$ to \eqref{main_res}, we have
\[ 0 \to \Hom_Q(S,S') \to I_\odd^* \otimes \Gamma(S') \xrightarrow{\phi^*} I_\ev^* \otimes \Gamma(S'(1)). \]
Taking global sections of \eqref{main_res} and its twist by $\OO(1)$, we can augment this to
\[ \xymatrix{
 & & & 0 \ar[d] \\
 & & & I_\ev^* \otimes I'_\ev \ar[d]^{\phi'_*} \\
 & & I_\odd^* \otimes I'_\odd \ar[r]^-{\phi^*} \ar@{=}[d] & I_\ev^* \otimes I'_\odd \otimes V^* \ar[d] \\
0 \ar[r] & \Hom_Q(S,S') \ar[r] & I_\odd^* \otimes \Gamma(S') \ar[r]^-{\phi^*} & I_\ev^* \otimes \Gamma(S'(1)) \ar[d] \\
 & & & 0
} \]
where the bottom row and the right column are exact.  Thus $\Hom_Q(S,S')$ is the set of $A \in \Hom_\CC(I_\odd,I'_\odd)$ for which there is a $B \in \Hom_\CC(I_\ev,I'_\ev)$ with $A \phi = \phi' B$; here we are thinking of $A$ and $B$ as matrices of complex numbers and $\phi$ and $\phi'$ as matrices of linear forms.  Since $\phi'_*$ is injective, such a $B$ is unique.  Multiplying $A \phi = \phi' B$ by $\psi$ on the right and $\psi'$ on the left, we have $\psi' A q = q B \psi$, so $\psi' A = B \psi$.  Thus
\[ \begin{pmatrix} B & 0 \\ 0 & A \end{pmatrix}
\begin{pmatrix} 0 & \psi \\ \phi & 0 \end{pmatrix}
= \begin{pmatrix} 0 & \psi' \\ \phi' & 0 \end{pmatrix}
\begin{pmatrix} B & 0 \\ 0 & A \end{pmatrix}. \]
Now $\left( \begin{smallmatrix} 0 & \psi \\ \phi & 0 \end{smallmatrix} \right)$ is a matrix of linear forms, and plugging in any $v \in V$ we get the map $I \to I$ given by left multiplication by $v$.  The $v$s generate $\Cl$, so $\left( \begin{smallmatrix} B & 0 \\ 0 & A \end{smallmatrix} \right) $ is in fact a homomorphism of $\Cl$-modules, not just of vector spaces.  Since the matrix is block diagonal, it respects the grading.
\end{proof}

%% file: dual.tex

\section{The Dual Sheaf}\label{dual}

To understand the dual sheaf $S^* := \HH om(S,\OO_Q)$, we work with some resolutions of $S$ on $Q$.  Restricting \eqref{main_res} to $Q$, we get
\[ \OO_Q^N(-1) \xrightarrow\phi \OO_Q^N \to S \to 0. \]
The sequence
\[ \dotsb \to \OO_Q^N(-2) \xrightarrow\psi \OO_Q^N(-1) \xrightarrow\phi \OO_Q^N \xrightarrow\psi \OO_Q^N(1) \xrightarrow\phi \OO_Q^N(2) \to \dotsb \]
on $Q$ is exact because $\phi \psi$ is a matrix factorization of $q$.  We can break it into exact sequences
\begin{multline*}
\dotsb \to \OO_Q^N(-3) \xrightarrow\phi \OO_Q^N(-2) \xrightarrow\psi \OO_Q^N(-1) \xrightarrow\phi \OO_Q^N \to S \to 0 \\
0 \to S \to \OO_Q^N(1) \xrightarrow\phi \OO_Q^N(2) \xrightarrow\psi \OO_Q^N(3) \xrightarrow\phi \OO_Q^N(4) \to \dotsb.
\end{multline*}
The former is a resolution of $S$ by $\HH om_Q(-, \OO_Q)$-acyclics, and applying \linebreak $\HH om_Q(-, \OO_Q)$ we get
\[ 0 \to S^* \to \OO_Q^N \xrightarrow{\phi^*} \OO_Q^N(1) \xrightarrow{\psi^*} \OO_Q^N(2) \xrightarrow{\phi^*} \OO_Q^N(3) \to \dotsb, \]
where $\phi^*$ and $\psi^*$ are the transposes of $\phi$ and $\psi$.  This is exact because $\psi^* \phi^*$ is also a matrix factorization of $q$.  Thus $S^*$ is the cokernel of
\begin{equation}\label{res_of_dual}
\OO_{\PP V}(-2) \otimes I_\odd^* \xrightarrow{\phi^*} \OO_{\PP V}(-1) \otimes I_\ev^*
\end{equation}
and $\EE xt^{>0}_Q(S, \OO_Q) = 0$.  That is, $S^*$ is not just the dual of $S$, but the derived dual.  Also, we observe that $S^{**} = S$; that is, $S$ is reflexive.

From all this we suspect that $S^*(1)$ is a spinor sheaf.  In fact it is:
\begin{proposition}
If $\codim W$ is odd then $S^* \cong S(-1)$ as $G_\ev$-equivariant sheaves.  If $\codim W$ is even then $S^* \cong T(-1)$.
\end{proposition}
\begin{proof}
Let ${}^\top$ be the anti-automorphism of $\Cl$ determined by $(v_1 \dotsm v_k)^\top = v_k \dotsm v_1$.  Then $I^\top$ is the right ideal $w_1 \dotsm w_m \cdot \Cl$.  The dual vector space $I^*$ is a right $\Cl$-module via the action $(f \cdot \xi)(-) = f(\xi-)$.  We will show that these two right modules are isomorphic up to a shift.

The natural filtration of the tensor algebra descends to $\Cl$,
\[ 0 = F_0 \subset F_1 \subset \dotsb \subset F_{\dim V} = \Cl, \]
and the associated graded pieces are $F_i/F_{i-1} = \textstyle\bigwedge^i V$.  In particular $\Cl/F_{\dim V-1}$ is 1-dimensional, so by choosing a generator we get a linear form $\tr: \Cl \to \CC$.  The pairing $\Cl \otimes \Cl \to \CC$ given by $\xi \otimes \eta \mapsto \tr(\xi\eta)$ is non-degenerate.  If $v \in V$ and $\xi \in \Cl$ then $\tr(v\xi) = \pm\tr(\xi v)$.

We claim that $I^*$ is generated by $\tr|_I$ and is isomorphic as an \emph{ungraded} module to $I^\top$.  Since $\dim I^* = \dim I^\top$, it is enough to check that $w_1 \dotsm w_m$ and $\tr|_I$ have the same annihilator.  If $\tr|_I \cdot \xi = 0$, that is, if $\tr(\xi \eta w_1 \dotsm w_m) = 0$ for all $\eta \in \Cl$, then $\tr(w_1 \dotsm w_m \xi \eta) = 0$ for all $\eta$, so $w_1 \dotsm w_m \xi = 0$; and conversely.

Now $\tr|_I$ has degree $\dim V \pmod2$, and $w_1 \dotsm w_m$ has degree $\dim W \linebreak \pmod2$, so $I^*$ is isomorphic to $I^\top$ or $I^\top[1]$ according as $\codim W$ is even or odd.  So if $\codim W$ is even then \eqref{res_of_dual} becomes
\[ \begin{array}{rclcrcl}
\OO_{\PP V}(-1) & \otimes & I_\odd^\top & \to & \OO_{\PP V} & \otimes & I_\ev^\top \\[2pt]
v & \otimes & \xi & \mapsto & 1 & \otimes & \xi v,
\end{array} \]
so from the isomorphism
\[ \xymatrix{
\OO_{\PP V}(-1) \otimes I_\odd \ar[r] \ar[d]_{1 \otimes {}^\top} &
\OO_{\PP V} \otimes I_\ev \ar[d]^{1 \otimes {}^\top} \\
\OO_{\PP V}(-1) \otimes I_\odd^\top \ar[r] & \OO_{\PP V} \otimes I_\ev^\top
} \]
we see that $S^*(1) \cong T$.  Similarly, if $\codim W$ is odd then $S^*(1) \cong S$.

These isomorphisms are equivariant, as follows.  In \eqref{equivariant}, an element $g \in G_\ev$ acts on $I$ by left multiplication by $g$, so it acts on $I^*$ by right multiplication by $g^{-1}$ and on $I^\top$ by right multiplication by $g^\top$.  But from the definition of $G$ we see that $g^\top g = 1$.
\end{proof}

%% file: horrocks.tex

\section{Linear Sections and Cones}\label{horrocks}

Any singular quadric can be described as a linear section of a smooth quadric or as a cone over a smooth quadric, so we would like to know what happens when spinor bundles are restricted to linear sections or pulled back via projection from the vertex of the cone.  We will see that the latter gives spinor sheaves corresponding to maximal linear spaces, while the former gives spinor sheaves corresponding to linear spaces that are maximal in the smooth locus of $Q$.  Thus our spinor sheaves corresponding to linear spaces that meet $Q_\text{sing}$ in interesting ways interpolate between the two obvious ways of extending spinor bundles to singular quadrics.

\begin{proposition}\label{linear_sections}
Suppose that $U \subset V$ is transverse to $W$.  Let $Q' = Q \cap \PP U$, and let $S'$ be the spinor sheaf on $Q'$ corresponding to $W \cap U \subset U$.  Then
\[ S' \cong \begin{cases}
S|_{Q'} & \text{if $\codim U$ is even,} \\
T|_{Q'} & \text{if $\codim U$ is odd.}
\end{cases} \]
\end{proposition}
\begin{proof}
Note that transversality is necessary to make $\rank S' = \rank S$.  Since we have some freedom to move $W$ without changing $S$, it is not a large restriction.

Let $I'$ be the $\Cl(U)$-module corresponding to $W \cap U$.  Let $I$ be the $\Cl(V)$-module corresponding to $W$; then $I$ is also a $\Cl(U)$-module since $\Cl(U)$ is a subring of $\Cl(V)$.  Restriction is right exact, so we have 
\[ \begin{array}{rclcrclcccc}
\OO_{\PP U}(-1) & \otimes & I_\ev & \to & \OO_{\PP U} & \otimes & I_\odd & \to & S' & \to & 0 \\[2pt]
u & \otimes & \xi & \mapsto & 1 & \otimes & f(u)\xi,
\end{array} \]
so it is enough to show that $I \cong I'[\codim U]$ as $\Cl(U)$-modules.

Choose a basis $u_1, \dotsc, u_l$ of $W \cap U$ and extend it to a basis $u_1, \dotsc u_l$, $w_{l+1}, \dotsc, w_m$ of $W$.  Then 
\begin{align*}
I' &= \Cl(U) \cdot u_1 \dotsm u_l, \\
I &= \Cl(V) \cdot u_1 \dotsm u_l \cdot w_{l+1} \dotsm w_m,
\end{align*}
and the map $I' \xrightarrow{\cdot w_{l+1} \dotsm w_m} I[\codim U]$ is an isomorphism.
\end{proof}

Together with \cite[Thm.~2.11]{ottaviani_spinor} this implies that when $q$ is non-degenerate and $W$ is maximal, our $S$ and $T$ are indeed the classical spinor bundles.

It also allows us to prove an analogue of Horrocks' splitting criterion.  Horrocks' criterion is a generalization of Grothendieck's theorem that every vector bundle on $\PP^1$ is a sum of line bundles.  It states that a vector bundle $E$ on $\PP^n$ is a sum of line bundles if and only if it is arithmetically Cohen-Macaulay (ACM), that is, $H^i(E(t)) = 0$ for $0 < i < n$ and all $t$; it is proved by induction on $n$.  Ottaviani \cite{ottaviani_horrocks} showed that a bundle $E$ on a smooth quadric 3-fold\footnote{On a quadric surface $Q$ we have to worry about whether we really mean ``sum of line bundles,'' since $\operatorname{Pic}(Q) \ne \operatorname{Pic}(\PP^3)$.} is a sum of line bundles if and only if $E$ and $E \otimes S$ are ACM, where $S$ is the spinor bundle, and observed that the same induction argument carries the result to higher-dimensional smooth quadrics.  Ballico \cite{ballico} observed that if a singular quadric $Q$ is a linear section of a larger smooth quadric and $S$ is the restriction of a spinor bundle, then the same induction argument still works.  Here we observe that it works with any of our spinor sheaves:

\begin{proposition}
Suppose that $\rank Q \ge 5$.  Let $S$ be any spinor sheaf on $Q$.  Then a vector bunle $E$ on $Q$ is a sum of the line bundles if and only if $E$ and $E \otimes S$ are ACM.
\end{proposition}
\begin{proof}
The ``only if'' statement follows from \eqref{main_res}.  For the ``if'' statement, we induct on $\dim Q$.  Suppose that $S$ corresponds to an isotropic subspace $W$.  If $\dim Q = 3$ then $Q$ is smooth and $S$ is either the classical spinor bundle (if $W$ is maximal) or the sum of several copies of it (if $W$ is smaller).

If $\dim Q > 3$, let $H$ be a hyperplane transverse to $\PP W$ and to $Q_\text{sing}$.  Let $Q' = Q \cap H$, whose rank is again at least 5.  From the exact sequence
\[ 0 \to E(-1) \to E \to E|_{Q'} \to 0 \]
we see that $E|_{Q'}$ is ACM, as is $(E \otimes S)|_{Q'} = E|_{Q'} \otimes S'$, where $S'$ is the spinor sheaf on $Q'$ corresponding to $\PP W \cap H$.  Thus $E|_{Q'}$ is a sum of line bundles $\bigoplus \OO_{Q'}(i)^{n_i}$.  Let $F = \bigoplus \OO_Q(i)^{n_i}$, and let $f: F|_{Q'} \to E|_{Q'}$ be an isomorphism.  Applying $\Hom(F, -)$ to the sequence above we get
\[ \Hom(F,E) \to \Hom(F, E|_{Q'}) \to \Ext^1(F,E(-1)). \]
The third term vanishes because $E$ is ACM and $F$ is a sum of line bundles, so we can extend $f$ to a map $\tilde f: F \to E$.  Now the line bundles $\det F|_{Q'}$ and $\det E|_{Q'}$ are isomorphic, so $\det F$ and $\det E$ are isomorphic, so $\det \tilde f$ is a non-zero section of $(\det F)^{-1} \otimes \det E \cong \OO_Q$ and hence does not vanish, so $\tilde f$ is an isomorphism.
\end{proof}

We finish the section with cones:

\begin{proposition}
Suppose that $U \subset W \cap K$.  Observe that $Q$ is a cone over a quadric $Q' \subset \PP(V/U)$, with vertex $\PP U$.  Let $p: Q \setminus \PP U \to Q'$ denote projection from the vertex, and $S'$ the spinor sheaf on $Q'$ corresponding to $W/U$.  Then
\[ p^* S' \cong \begin{cases}
S|_{Q \setminus \PP U} & \text{if $\dim U$ is even,} \\
T|_{Q \setminus \PP U} & \text{if $\dim U$ is odd.}
\end{cases} \]
\end{proposition}
\begin{proof}
Let $I'$ be the $\Cl(V/U)$-module corresponding to $W/U$, which is also a $\Cl(V)$-module via the map $\Cl(V) \to \Cl(V/U)$.  As usual let $I$ be the $\Cl(V)$-module corresponding to $W$.  Since $p^*$ is right exact, as in the proof of Proposition \ref{linear_sections} it is enough to show that $I' \cong I[\dim U]$ as a $\Cl(V)$-modules.

Choose a section $s: V/U \to V$ of the projection $V \to V/U$.  Then $s$ respects the quadratic forms, hence makes $\Cl(V/U)$ a subring of $\Cl(V)$, and as in the proof of Proposition \ref{linear_sections} we get an isomorphism $I \cong I'[\dim U]$ of $\Cl(V/U)$-modules, hence of $\Cl(V)$-modules.
\end{proof}

%% file: stable.tex

\section{Stability}\label{stable}

In this last section we show that $S$ and $T$ are stable if $W$ is maximal and properly semi-stable otherwise.  For background on stability we refer to Huybrechts and Lehn \cite{huybrechts_lehn}.  Recall that the \emph{slope} of a torsion-free sheaf $\EE$ on polarized variety $X$ is
\[ \mu(\EE) = \frac{\deg\EE}{\rank\EE}. \]
It can be read from the Hilbert polynomial of $\EE$: if
\[ \chi(\OO_X(t)) = \deg X \cdot \frac{t^n}{n!} + C \cdot \frac{t^{n-1}}{(n-1)!} + \dotsb \]
then
\[ \chi(\EE(t)) = \deg X \cdot \rank\EE \cdot \frac{t^n}{n!} + (C \cdot \rank\EE + \deg\EE) \cdot \frac{t^{n-1}}{(n-1)!} + \dotsb. \]
For example, from \eqref{main_res} we find that the slope of any spinor sheaf is 1.

If $\FF \subset \EE$ is a proper, saturated subsheaf (that is, $\EE/\FF$ is torsion-free) then either 
\begin{align*}
\mu(\FF) &< \mu(\EE) < \mu(\EE/\FF) \text{ or} \\
\mu(\FF) &= \mu(\EE) = \mu(\EE/\FF) \text{ or} \\
\mu(\FF) &> \mu(\EE) > \mu(\EE/\FF).
\end{align*}
We say that $\EE$ is \emph{semi-stable} if $\mu(\FF) \le \mu(\EE)$ for all such $\FF$ and \emph{stable} if $\mu(\FF) < \mu(\EE)$.  Every sheaf has a unique \emph{Harder-Narasimhan filtration}
\[ 0 = \FF_0 \subset \FF_1 \subset \dotsb \subset \FF_n = \EE, \]
where the quotients $\FF_i / \FF_{i-1}$ are semi-stable and $\mu(\FF_i/\FF_{i-1}) > \mu(\FF_{i+1}/\FF_i)$.  Every semi-stable sheaf has a \emph{Jordan-H\"older filtration}
\[ 0 = \FF_0 \subset \FF_1 \subset \dotsb \subset \FF_n = \EE, \]
where the quotients $\FF_i / \FF_{i-1}$ are stable and $\mu(\FF_i/\FF_{i-1}) = \mu(\EE)$.  This is not unique, but the associated graded object $\bigoplus_{i=1}^n \FF_i/\FF_{i-1}$ is.  Two semi-stable sheaves whose Jordan-H\"older filtrations have the same associated graded object are called \emph{S-equivalent}.  A sheaf is called \emph{polystable} if it is a direct sum of stable sheaves of the same slope.  \vspace\baselineskip

If $W$ is maximal then we will show in a moment that $S$ and $T$ are stable.  If $W'$ is codimension 1 in $W$ then \eqref{jordan_hoelder} gives a Jordan--H\"older filtration $0 \subset S \subset S'$ with $S/0 = S$ and $S'/S = T$, so $S'$ is properly semi-stable and S-equivalent to the polystable sheaf $S \oplus T$, as is $T'$.  If $W''$ is codimension 1 in $W'$ then the Jordan--H\"older filtration is slightly more complicated, but $S''$ and $T''$ are S-equivalent to $S \oplus S \oplus T \oplus T$.  In general the S-equivalence class of a spinor sheaf depends only on the dimension of the isotropic space.

Our proof that $S$ and $T$ are stable when $W$ is maximal will use the fact that they are simple.  We begin with a lemma:
\begin{lemma}
If $W$ is maximal then $I$ is irreducible.
\end{lemma}
\begin{proof}
If $\dim V/K = 2k$ is even, there is a basis $v_1, \dotsc, v_n$ of $V$ in which
\[ q = x_1 x_{k+1} + \dotsb + x_k x_{2k} \]
and $W = \operatorname{span}(v_{k+1}, \dotsc, v_n)$.  Let $\xi = v_{k+1} \dotsm v_n$ be the generator of $I$.  Then any $\xi' \in I$ different from zero is of the form
\[ \xi' = \alpha v_{i_1} \dotsm v_{i_l} \xi + \text{terms of the same or shorter length}, \]
where $\alpha \in \CC$ is not zero and $1 \le i_1 < \dotsb < i_l \le k$.  I claim that $v_{i_l+k} \dotsm v_{i_1+k} \xi' \linebreak = \alpha \xi$.  To see this, observe that if $1 \le i, j \le k$ then $v_i$ anti-commutes with $v_{j+k}$ when $i \ne j$ and that $v_{i+k} \xi = 0$, so left multiplication by $v_{i_l+k} \dotsm v_{i_1+k}$ annihilates any basis vector not containing $v_{i_1} \dotsm v_{i_l}$; and $v_{i+k} v_i \xi = (1 - v_i v_{i+k}) \xi = \xi$.  Thus any non-zero element of $I$ generates $I$, so $I$ is irreducible.

If $\dim V/K = 2k+1$ is odd, there is a basis $v_0, \dotsc, v_n$ of $V$ in which
\[ q = x_0^2 + x_1 x_{k+1} + \dotsb + x_k x_{2k} \]
and $W = \operatorname{span}(v_{k+1}, \dotsc, v_n)$.  Let $\xi = v_{k+1} \dotsm v_n$ be the generator of $I$.  Let $J \subseteq I$ be a graded submodule.  By an argument similar to the one given above, for any non-zero $\xi' \in J$ there are $1 \le i_1 < \dotsb < i_l \le k$ such that $v_{i_l} \dotsm v_{i_1} \xi' = (\alpha + \beta v_0) \xi$, where $\alpha, \beta \in \CC$ are not both zero.  Since $J$ is graded, it contains both $\alpha \xi$ and $\beta v_0 \xi$.  If $\alpha \ne 0$ then $\xi \in J$, and if $\beta \ne 0$ then $v_0 \cdot \beta v_0 \xi = \beta \xi$, so again $\xi \in J$, so $J = I$.
\end{proof}

\begin{proposition}
If $W$ is maximal then $S$ is simple, that is, $\Hom_Q(S,S) = \CC$.  If $\dim V/K$ is even, $\pi(W)$ is maximal in $V/K$, and $W \cap K$ is codimension 1 in $K$, then again $S$ is simple.  Otherwise $S$ is not simple.
\end{proposition}
\begin{proof}
The first statement is immediate from the previous lemma, Schur's \linebreak lemma, and Proposition \ref{fully_faithful}.

For the second statement, let $W' = W+K$.  Let $J$ be a proper submodule of $I$, and consider the short exact sequence
\[ 0 \to I' \to I \xrightarrow{\cdot v} I'[1] \to 0, \]
where $v \in K \setminus W$.
Since $I'$ is irreducible, either $J \cap I' = 0$ or $J \supset I'$.  If $J \cap I' = 0$ then $Jv$ is isomorphic to $J$; since $I'[1]$ is irreducible, either $Jv = 0$, so $J = 0$, or $Jv = I'[1]$, so $I = I' \oplus J = I' \oplus I'[1]$, which we know is not true.  If $J \supset I'$ then again either $Jv = 0$, so $J = I'$, or $Jv = I'[1]$, so $J = I$.  Thus the only proper submodule of $I$ is $I'$, and the only proper quotient is $I'[1]$.  Since these are not isomorphic, any homomorphism $I \to I$ is an isomorphism or zero, so again by Schur's lemma $\Hom_\Cl(I,I) = \CC$.

For the third statement, if $\pi(W)$ is not maximal in $V/K$ then $S$ is a direct sum, hence is not simple.  If $W \cap K$ is codimension 2 or more in $K$, choose $W'' \supset W' \supset W$ with $\pi(W'') = \pi(W') = \pi(W)$; then the composition
\[ S \twoheadrightarrow T' \twoheadrightarrow S'' \hookrightarrow S' \hookrightarrow S \]
is neither zero nor an isomorphism.  If $W \cap K$ is codimension 1 in $K$ and $\dim V/K$ is odd, let $W' = W + K$; then the composition
\[ S \twoheadrightarrow T' \cong S' \hookrightarrow S \]
is neither zero nor an isomorphism.
\end{proof}

\begin{theorem*}
Suppose that $\rank Q > 2$.  If $W$ is maximal then $S$ and $T$ stable.
\end{theorem*}
\begin{proof}
It is enough to show this for $S$.  Suppose that a subsheaf $\FF \subset S$ is invariant under the action of $G_\ev$ introduced in \S\ref{dependence}.  If $\rank Q > 2$ then $G_\ev$ acts transitively on the smooth locus of $Q$, so $\FF$ is a vector bundle there.  Let $p \in Q$ be a smooth point and $H \subset G_\ev$ be the stabilizer of $p$; then according to Ottaviani \cite{ottaviani_spinor}, the representation of $H$ on the fiber $S|_p$ is irreducible (recall that $G_\ev$ contains $\Spin(U, q|_U)$ for any $U$ complementary to $K$).  Thus either $\rank \FF = 0$, so $\FF = 0$ since $S$ is reflexive and hence torsion-free, or $\rank \FF = \rank S$.

Thus $S$ has no invariant proper saturated subsheaves.  The Harder--\linebreak Narasimhan filtration is unique, hence invariant, so $S$ is semi-stable.  Consider the \emph{socle} of $S$, that is, the maximal polystable subsheaf, which is necessarily saturated.  This too is unique, hence invariant, hence is $S$; that is, $S$ is a direct sum of stable sheaves.  But $S$ is simple, hence indecomposable, so $S$ is stable.
\end{proof}

If $\rank Q = 2$ then $Q$ is a union of hyperplanes $H$ and $H'$ and $S$ is either $\OO_H$ or $\OO_{H'}$, which are torsion and thus not eligible for slope stability, but they have no proper saturated subsheaves and thus are Gieseker stable.  We excluded from the beginning the non-reduced case $\rank Q = 1$.